\documentclass{amsart}
\usepackage{amssymb}
\usepackage{amscd,enumerate,amsfonts,calc,amsmath,verbatim}

\newcommand{\mlabel}{\label}

\newcommand{\bd}{\boldsymbol}
\newcommand{\fm}{\mathfrak{m}}
\newcommand{\fn}{\mathfrak{n}}
\newcommand{\fp}{\mathfrak{p}}
\newcommand{\fs}{\mathfrak{s}}
\newcommand{\ft}{\mathfrak{t}}
\newcommand{\ov}{\overline}
\newcommand{\wh}{\widehat}

\newcommand{\codim}{\operatorname{codim}}
\newcommand{\depth}{\operatorname{depth}}
\newcommand{\edim}{\operatorname{edim}}
\newcommand{\mult}{\operatorname{mult}}
\newcommand{\rank}{\operatorname{rank}}
\newcommand{\supp}{\operatorname{Supp}}

\newcommand{\pd}{\operatorname{pd}}
\newcommand{\id}{\operatorname{id}}

\newcommand{\Hom}{\operatorname{Hom}}
\newcommand{\Ext}{\operatorname{Ext}}
\newcommand{\Tor}{\operatorname{Tor}}
\newcommand{\Po}{\operatorname{P}}
\newcommand{\Ba}{\operatorname{I}}
\newcommand{\Hi}{\operatorname{H}}
\newcommand{\CIdim}{\operatorname{CI-dim}}
\newcommand{\Syz}{\operatorname{Syz}}
\newcommand{\V}{\operatorname{V}}
\newcommand{\GTE}{\operatorname{\mathbf{GTE}}}
\newcommand{\G}{\operatorname{\mathbf G}}

\newcommand{\GGO}{\operatorname{\mathbf {GGO}}}
\newcommand{\GH}{\operatorname{\mathbf {GH}}}

\swapnumbers

\theoremstyle{plain}
\newtheorem{theorem}{Theorem}[section]

\newtheorem*{TheoremAB}{Theorem AB}
\newtheorem*{TheoremHJ}{Theorem HJ}

\newtheorem*{Theorem1}{Theorem 1}
\newtheorem*{Theorem3}{Theorem 3}
\newtheorem*{Theorem2}{Theorem 2}

\newtheorem{proposition}[theorem]{Proposition}

\newtheorem{lemma}[theorem]{Lemma}

\newtheorem{corollary}[theorem]{Corollary}

\theoremstyle{definition}

\newtheorem{definition}[theorem]{Definition}
\newtheorem{remark}[theorem]{Remark}

\newtheorem{example}[theorem]{Example}

\newtheorem{chunk}[theorem]{}
\newtheorem{subchunk}{}
\newtheorem*{chunk*}{}

\theoremstyle{remark}

\newtheorem*{caseone}{Case {\rm 1}}
\newtheorem*{casetwo}{Case {\rm 2}}

\numberwithin{equation}{theorem}
\numberwithin{subchunk}{theorem}

\begin{document}

\title[Vanishing of cohomology] {Vanishing of
 cohomology over\\ Gorenstein rings of small codimension}

\author[L.~M.~\c Sega]{Liana M.~\c Sega}

\address{Department of Mathematics, Purdue University, West Lafayette,
Indiana 47907}

\email{lmsega@math.purdue.edu} 

\date{}

\commby{Wolmer Vasconcelos}


\begin{abstract}

  We prove that if $M$, $N$ are finite modules over a Gorenstein local
  ring $R$ of codimension at most $4$, then the vanishing of
  $\Ext^n_R(M,N)$ for $n\gg 0$ is equivalent to the vanishing of
  $\Ext^n_R(N,M)$ for $n\gg 0$. Furthermore, if $\wh R$ has no
  embedded deformation, then such vanishing occurs if and only if $M$
  or $N$ has finite projective dimension.

\end{abstract}
\maketitle

 \section*{Introduction}

 Let $(R,\fm,k)$ be a commutative noetherian local ring with maximal
 ideal $\fm$ and residue field $k$. The {\em codimension} of $R$ is
 the number $\codim R=\edim R-\dim R$, where $\edim R$ denotes the
 minimal number of generators of $\fm$.

 Let $M$, $N$ be finite $R$-modules. We study the vanishing of
 $\Ext_R^n(M,N)$ for all $n\gg 0$ and the vanishing of $\Tor^R_n(M,N)$
 for all $n\gg 0$, assuming that $R$ is Gorenstein and $M$, $N$ are
 finite $R$-modules.

It is convenient to state the vanishing properties of (co)homology
in terms of numbers $e_R(M,N)$ and $t^R(M,N)$, defined as follows:
 \begin{align*}
 e_R(M,N)&=\sup\{n\in \mathbb N\mid \Ext_R^n(M,N)\ne 0 \}\\
 t^R\,(M,N)&=\sup\{n\in \mathbb N\mid \Tor^R_n(M,N)\ne 0 \}
 \end{align*}

A local ring $(Q,\fn,k)$ is called an {\em {\rm(}embedded{\rm)\/}
deformation} of $R$ if $R\cong Q/(\bd f)$, where $\bd f$ is a
$Q$-regular sequence (contained in $\fn^2$). The ring $R$ is said
to be a {\em complete intersection} if its $\fm$-adic completion
$\wh R$ has an embedded deformation which is a regular ring.
Vanishing of (co)homology over complete intersections has
remarkable properties, as proved by Avramov and Buchweitz
\cite[Theorem III]{AB}. They relate vanishing properties to
homogeneous algebraic varieties $V^*_R(M)\subseteq
\widetilde{k}^c$, defined in \cite{A2}, where $c=\codim R$ and
$\widetilde{k}$ denotes the algebraic closure of $k$.

 \begin{TheoremAB}
Let $R$ be a complete intersection local ring. One of the numbers
$e_R(M,N)$, $e_R(N,M)$, $t^R(M,N)$ is finite if and only if all of
them are less than or equal to $\dim R$, if and only if
$\V^*_R(M)\cap \V^*_R(N)=0$.
 \end{TheoremAB}

 An interesting consequence of Theorem AB is that the vanishing of
 $\Ext^n_R(M,N)$ for all $n\gg 0$ is symmetric in the two module
 variables. It is clear that a ring with this property is Gorenstein.
 To study the symmetry in the vanishing of $\Ext$, Huneke and
 Jorgensen \cite{HJ} define a class of Gorenstein rings, called {\em AB
   rings}, by the condition that $\sup\{e_R(M,N)\}$ is finite
 when the supremum is taken over all finite $R$-modules $M$, $N$ with
 $e_R(M,N)<\infty$. In \cite[3.1, 3.3(1), 4.1]{HJ} they prove:

\begin{TheoremHJ}
The ring $R$ is AB if and only if every pair $(M,N)$ of finite
$R$-modules with $e_R(M,N)<\infty$ {\rm (}respectively
$t^R(M,N)<\infty${\rm)} satisfies $e_R(M,N)\le \dim R$ and
$e_R(N,M)\le \dim R$ {\rm (}respectively $t^R(M,N)\le \dim R${\rm
)}.
\end{TheoremHJ}

Theorem AB shows that a complete intersection local ring is an AB
ring. A non-complete intersection example of AB rings is given by
Gorenstein rings of minimal multiplicity and codimension at least
$3$: It is proved in \cite[3.5]{HJ} that over such rings the
finiteness of $e_R(M,N)$ or $t^R(M,N)$ implies one of the modules
has finite projective dimension.

In Theorem 1 we identify another class of rings with this
property.

 \begin{Theorem1}
Let $R$ be a local Gorenstein ring such that $\codim R\le 4$ and
$\wh R$ admits no embedded deformation. If one of the numbers
$e_R(M,N)$, $t^R(M,N)$ is finite, then $M$ or $N$ has finite
projective dimension.
 \end{Theorem1}

The theorem is proved in Section 2 as Theorem \ref{T1}. The proof
uses the fact that, for the rings in Theorem 1, all finite
$R$-modules have rational Poincar\'e series which share an
explicitly known denominator: this is due to Jacobsson, Kustin,
Miller \cite{J}, \cite{JKM}, \cite{K}.  Preliminaries on rational
Poincar\'e series are presented in Section 1, where we also give a
new proof of the result on Gorenstein rings of minimal
multiplicity. Clearly, such an approach is not applicable to rings
that have finite modules with irrational Poincar\'e series. Such
rings have been constructed by B\o gvad \cite{B} in codimension
$12$ and higher.

To describe the aymptotic vanishing of cohomology over all
Gorenstein rings with $\codim R\le 4$, we use the notion of finite
complete intersection dimension introduced by Avramov, Gasharov
and Peeva \cite{AGP}. An $R$-module $M$ is said to have {\em
finite CI-dimension} if there exist a flat ring homomorphism $R\to
R'$ and a deformation $Q$ of $R'$ such that
$\pd_Q(M\otimes_RR')<\infty$. Note that any module of finite
projective dimension has finite CI-dimension, and that if $R$ is a
complete intersection, then every finite $R$-module has finite
CI-dimension. The next statement is derived from the theorems
presented above, as well as from expressions for $e_R(M,N)$ and
$t^R(M,N)$ obtained by Araya, Yoshino, Avramov, Buchweitz
\cite{AY}, \cite{AB}. It shows, in particular, that all Gorenstein
rings of codimension at most $4$ are AB:

\begin{Theorem2}
Let $R$ be a local Gorenstein ring with $\codim R\le 4$.

If $e_R(M,N)<\infty$, then $M$ or $N$ has finite CI-dimension and
\begin{align*}
e_R(M,N)&=\dim R-\depth_RM\\
e_R(N,M)&=\dim R-\depth_RN
\end{align*}

If $t^R(M,N)<\infty$, then $M$ or $N$ has finite CI-dimension and
\[
t^R(M,N)\le\max\{\dim R-\depth_RM,\dim R-\depth_RN\}\,.
\]
\end{Theorem2}

The theorem is proved as Theorem \ref{T2} in Section 3. An example
is provided there to show that the conclusion regarding
CI-dimension cannot be extended to high codimensions. However, it
is not known (to the author) whether the other conclusions can be
extended or not.

Under the hypotheses of Theorem 1, asymptotic vanishing of $\Tor$
or $\Ext$ occurs only for trivial reasons: one of the modules has
finite projective dimension. The situation is different if $R$ is
a complete intersection with $\codim R\ge 2$. Over such rings,
Avramov \cite[6.5]{A2} constructs for every linear subspace
$H\subseteq {\widetilde k}^c$ an $R$-module $M$ of finite length
with $\V^*_R(M)=H$; taking modules $M$, $N$ that correspond to
nontrivial subspaces with zero intersection, one sees from Theorem
AB that the numbers $e_R(M,N)$, $e_R(N,M)$ and $t^R(M,N)$ are all
finite.

 We show that (co)homology may vanish for nontrivial reasons whenever
 $\wh R$ has a non-regular embedded deformation:

 \begin{Theorem3}
   Let $R$ be a local Gorenstein ring with $\codim R\ge 2$. If
   $\wh R$ has an embedded deformation, then there exist $R$-modules
   $M$, $N$ of finite length with $\pd_RM=\pd_RN=\infty$ and
   $e_R(M,N)=e_R(N,M)=t_R(M,N)=\dim R$.
 \end{Theorem3}

This result is contained in the more general Corollary \ref{last}
proved in Section 4. The hypothesis on the codimension is
necessary. Indeed, if $R$ is a Gorenstein ring with $\codim R\le
1$, then $\wh R=Q/(f)$ for a regular local ring $Q$ and a
$Q$-regular element $f$. For such rings, $t^R(M,N)<\infty$ or
$e_R(M,N)<\infty$ implies that one of the modules $M$ or $N$ has
finite projective dimension, cf.\ Huneke and Wiegand
\cite[1.9]{HW}, respectively Avramov and Buchweitz
\cite[5.12]{AB}.

 \section{Rational Poincar\'e series}

 In this section $(R,\fm,k)$ is a local ring and $M$, $N$ denote
 finite $R$-modules. We use several formal power series with integer
 coefficients associated to a finite $R$-module.

 We let $\Po^R_M(t)$ denote the {\em Poincar\'e series\/} of $M$ over
 $R$
 \[
 \Po^R_M(t)=\sum_{i=0}^{\infty}\rank_k\Tor^R_i(M,k)t^i\in\mathbb Z[[t]]
 \]
 and $\Ba^M_R(M)$ denote the {\em Bass series\/} of $M$ over $R$
\[
 \Ba^M_R(t)=\sum_{i=0}^{\infty}\rank_k\Ext_R^i(k,M)t^i\in\mathbb Z[[t]]\,.\
\]
 We also use the {\em Hilbert series} of $M$ over $R$, defined to be
 the series
 \[
 \Hi_{M}^R(t)=\sum_{i=0}^{\infty}\rank_k(\fm^iM/\fm^{i+1}M)t^i\in\mathbb Z[[t]]\,.
 \]
 Recall that $\Hi_{M}^R(t)$ represents a rational function with
 denominator $(1-t)^{\dim M}$.

 \begin{chunk}
   \mlabel{completion} Let $R\to R'$ be a flat homomorphism of local
   rings such that $\fm R'$ is the maximal ideal of $R'$. For
   $M'=M\otimes_RR'$ and $N'=N\otimes_R R'$ there are equalities
   $\Po^{R}_{M}(t)=\Po^{R'}_{M'}(t)$ and
   $\Ba_{R}^{M}(t)=\Ba_{R'}^{M'}(t)$, $e_R(M,N)=e_{R'}(M' ,N')$ and
   $t^R(M,N)=t^{R'}(M',N')$, $\pd_{R'}(M')=\pd_RM$ and
   $\id_{R'}(N')=\id_RN$. This situation applies in particular to the
   completion map $R\to \wh R$.
 \end{chunk}

 \begin{lemma}
   \mlabel{Foxby} Assume there exists a polynomial $c(t)$ such that
   $c(t)\Po^R_M(t)\in \mathbb Z[t]$ for each finite $R$-module $M$. If
   $c(t)\Po^R_k(t)=(1+t)^ea(t)$ for an integer $e\ge 0$ and a
   polynomial $a(t)\in\mathbb Z[t]$ with $a(-1)\ne 0$, then
   $c(t)(1+t)^m\Ba_R^N(t)\in\mathbb Z[t]$ for each finite $R$-module
   $N$, where $m=\max\{0,\dim R-e\}$.
 \end{lemma}

 \begin{proof}
   Set $d=\dim R$. By Foxby \cite[3.10]{F} there exists a finite $\wh
   R$-module $L$ and a polynomial $g(t)$ with integer coefficients and
   degree less than $d$ such that
   \[
   \Ba_R^N(t)=g(t)+t^d\Po^{\wh R}_{L}(t)\,.
   \]
   As shown by Levin \cite[p.\,8]{L} there exists an integer $n$ such that
   $M=L/\fm^nL$ satisfies
   \[
   \Po^{\wh R}_{L}(t)=\Po^{\wh R}_{M}(t)-t\Hi_{\fm^nL}^{\wh
     R}(-t)\Po^{\wh R}_k(t)\,.
   \]
   The $\wh R$-module $M$ has finite length, so considered as an
   $R$-module it has finite length and is complete. From
   \ref{completion} we get $\Po^{R}_{M}(t)=\Po^{\wh R}_{M}(t)$ and
   $\Po^{R}_{k}(t)=\Po^{\wh R}_{k}(t)$.  Since $\dim_{\wh
     R}(\fm^nL)\le d$, we have $(1-t)^d\Hi_{\fm^nL}^{\wh
     R}(t)\in\mathbb Z[t]$ and the conclusion follows.
 \end{proof}

 The next result connects vanishing of (co)homology to Bass series and
 Poincar\'e series. Part (2) is due to C.  Miller, cf.\ the proof of
 \cite[1.1]{M}; part (1) is \cite[1.5(1)]{AB}.

 \begin{chunk}
   \mlabel{PB} For finite $R$-modules $M$, $N$ the following hold:

 \vspace{0.15cm}

 (1) If $e_R(M,N)=0$, then
 $\Ba_R^{\Hom_R(M,N)}(t)=\Po^R_M(t)\Ba_R^N(t)$

 \vspace{0.15cm}

 (2) If $\,t^R(M,N)=0$, then
 $\Po^R_{M\otimes_RN}(t)=\Po^R_M(t)\Po^R_N(t)$.
 \end{chunk}

 \begin{definition}
   \mlabel{prop} A factorization $c(t)=p(t)\cdot q(t)\cdot r(t)$ in
   $\mathbb Z[t]$ is said to be {\em good \/} if $p(t)=1$ or $p(t)$ is
   irreducible, $q(t)$ has non-negative coefficients, and $r(t)=1$ or
   $r(t)$ is irreducible and has no positive real root among its
   complex roots of minimal absolute value.
 \end{definition}

 \begin{proposition}
   \mlabel{basic} Let $R$ be a local ring for which there exists a
   polynomial $c(t)$ such that $c(t)\Po^R_L(t)\in\mathbb Z[t]$ for
   each finite $R$-module $L$. When $c(t)$ has a good factorization
   the following hold for all finite $R$-modules $M$, $N$:
 \begin{enumerate}[{\quad\rm(1)}]
 \item If $e_R(M,N)<\infty$, then $\pd_RM<\infty$ or
   $\id_RN<\infty$.
 \item If $t^R(M,N)<\infty$, then $\pd_RM<\infty$ or
   $\pd_RN<\infty$.
 \end{enumerate}
 \end{proposition}

 \begin{proof}
   We only give a proof of (1); the proof of (2) is similar.

   Replacing $M$ by a syzygy, if necessary, we may assume
   $e_R(M,N)=0$. Set $b(t)=(1+t)^dc(t)$, where $d=\dim R$. By
   hypothesis and Lemma \ref{Foxby} we have
\[
b(t)\Po^R_M(t)=m(t)\in\mathbb Z[t]\quad\text{and}\quad b(t)\Ba_R^N(t)=n(t)\in\mathbb Z[t]\,.
 \]
 Since $(1+t)^d$ is in $\mathbb N[t]$ it follows that $b(t)$ has a
 good factorization $b(t)=p(t)\cdot q(t)\cdot r(t)$ as in Definition
 \ref{prop}.  Using \ref{PB}(1) we obtain:
   \[
   b(t)\Ba_R^{\Hom_R(M,N)}(t)=b(t)\Po^R_M(t)\Ba_R^N(t)=\frac{m(t)n(t)}{p(t)q(t)r(t)}\,.
   \]
   By Lemma \ref{Foxby} the expression on the left is a polynomial, so
   $p(t)$ divides $m(t)n(t)$. Since $p(t)=1$ or $p(t)$ is irreducible,
   two cases arise.

   \caseone $m(t)=p(t)m_1(t)$ with $m_1(t)\in \mathbb Z[t]$. We then have an equality
\[
q(t)\Po^R_M(t)=\frac{m_1(t)}{r(t)}\,.
\]
Assume $r(t)$ does not divide $m_1(t)$. This means that the radius of
convergence $\rho$ of these power series is finite. The series
$q(t)\Po^R_M(t)$ has non-negative coefficients, so by the Pringsheim Principle, cf.\ \cite[7.2]{Ti}, $\rho$ is a singular point for it. On the
other hand, the hypothesis on $r(t)$ implies that $\rho$ is not a
singular point for $m_1(t)/r(t)$. This contradiction shows that $r(t)$
divides $m_1(t)$, so $q(t)\Po^R_M(t)$ is a polynomial.  As both
factors have non-negative coefficients, we conclude that $\Po^R_M(t)$
is a polynomial, that is, $\pd_R(M)<\infty$.

\casetwo $n(t)=p(t)n_1(t)$ with $n_1(t)\in \mathbb Z[t]$. An argument
similar to the one above shows that $\Ba_R^N(t)$ is a polynomial,
hence $\id_RN<\infty$.
\end{proof}

The two parts of the proposition can be combined when the ring is Gorenstein:

\begin{chunk}
  \mlabel{Gor} If $R$ is Gorenstein, then $\id_RN<\infty$ is
  equivalent to $\pd_RN<\infty$.
\end{chunk}

Let $R$ be a Gorenstein ring with $\codim R\ge 2$. It is known that
$R$ has multiplicity at least $\codim R+2$, cf.\ \cite[3.2]{Sa}. When
equality holds, $R$ is said to be {\em Gorenstein of minimal
  multiplicity}.

\begin{remark}
  \mlabel{multiplicity} If $R$ is Gorenstein of minimal
  multiplicity with $\codim R\ge 2$, then $c(t)\Po^R_M(t)\in \mathbb Z[t]$ for each
  finite $R$-module $M$, where $c(t)=1-(\codim R)t+t^2$.

  To see this, we may assume that $k$ is infinite: if not, then
  $R'=R[t]_{\fm[t]}$ has $\mult R'=\mult R$ and we apply
  \ref{completion}. It also suffices to consider maximal
  Cohen-Macaulay modules $M$. If $\dim R=0$, then $\fm^3=0$ and the
  result is proved by Sj\"odin \cite{Sj}. If $\dim R>0$, then there
  exists a maximal $R$-regular sequence $\bd g$ such that
  $\codim R=\codim R/(\bd g)$ and $\mult R=\mult R/(\bd g)$. Since $M$ is
  maximal Cohen-Macaulay, this sequence is also $M$-regular.
  The isomorphisms $\Tor^R_n(M,k)\cong\Tor^{R/(\bd g)}_n(M/\bd gM,k)$
  then yield $\Po ^R_M(t)= \Po^{R/(\bd g)}_{M/\bd gM}(t)$. The
  conclusion follows by the dimension $0$ case.
 \end{remark}

 If $\codim R>2$, then the polynomial $c(t)=1-(\codim R)t+t^2$ is
 irreducible, hence $c(t)=c(t)\cdot 1\cdot 1$ is a good factorization.
 From Proposition \ref{basic} and Remark \ref{multiplicity}, we obtain
 a new proof of \cite[3.5]{HJ}:

\begin{corollary}
  \mlabel{example} Let $R$ be a local Gorenstein ring of minimal
  multiplicity with $\codim R\ge 3$. If $M$, $N$ are finite
  $R$-modules such that one of the numbers $e_R(M,N)$, $t^R(M,N)$ is
  finite, then $M$ or $N$ has finite projective dimension.\qed
\end{corollary}

 \section{Gorenstein rings of small codimension}

 In this section $R$ denotes a Gorenstein ring of codimension at most
 $4$ which is not a complete intersection.  Over such rings, the
 Poincar\'e series of all finite $R$-modules are rational and share a
 common denominator. This result, and the form of the denominator were
 obtained by Jacobsson \cite{J} in codimension $3$, by Jacobsson,
 Kustin and Miller \cite[2.3]{JKM} in codimension $4$, characteristic
 different from $2$ and by Kustin \cite{K} in codimension $4$ and
 characteristic different from $3$. We collect below the relevant information;
 we refer to \cite[\S 3]{A1} for details.

\begin{chunk}
  \mlabel{classify} There exists a polynomial $c(t)\in\mathbb Z[t]$
  such that $c(t)\Po^R_M(t)\in\mathbb Z[t]$ for all finite $R$-modules
  $M$.  This polynomial has the form $c(t)=d(t)(1+t)^m$, where the
  polynomial $d(t)\in\mathbb Z[t]$ and the non-negative integer
  $m$ are as follows:

 \begin{center}
\renewcommand{\arraystretch}{1.5}
  \begin{tabular}{|c|c|c|c|c|}
    \hline
    type       &$\codim R$ &$d(t)$             &$m$               & restrictions\\ \hline
    $\G(l+1)$  &3 &$1-t-lt^2-t^3+t^4$        &1       & $l\ge 4$ \\ \hline
    $\GTE$    &4 & $1-2t-(l-2)t^2+t^3+t^4-t^5$  &2     & $l\ge 5$ \\ \hline
    $\GGO$   &4 & $1-2t-(l-2)t^2-2t^3+t^4$      &2    & $l\ge 5$\\ \hline
    $\GH(p)$ &4 & $1-2t-(l-2)t^2+(p-2)t^3+2t^4-t^5$ &2 & $1\le p\le l\ge 5$ \\ \hline
  \end{tabular}
\end{center}
\end{chunk}

Avramov \cite[3.1]{A1} determines when $\wh R$ has an embedded
deformation:

\begin{chunk}
  \mlabel{embedded} The ring $\wh R$ admits an embedded deformation if
  and only if $d(1)=0$, if and only if $R$ is of type $\GH(p)$,
  with $p=l$.
 \end{chunk}

 The proof of the next result goes through a careful
 examination of the polynomials in the chart above. Earlier, Sun
 \cite{Su} used a case by case analysis of these polynomials to prove
 that the Betti numbers of finite modules over a Gorenstein ring of
 codimension at most $4$ are eventually non-decreasing.

\begin{theorem}
\mlabel{T1}
   Let $R$ be a local Gorenstein ring such that $\codim R\le 4$ and
   $\wh R$ admits no embedded deformation. If $M$, $N$ are finite
   $R$-modules and one of the numbers $e_R(M,N)$, $t^R(M,N)$ is
   finite, then $M$ or $N$ has finite projective dimension.
 \end{theorem}

 \begin{proof}
   By Proposition \ref{basic} and \ref{Gor}, it suffices to show that the common
   denominator $c(t)$ has a good factorization. Since $(1+t)^m$ has
   non-negative coefficients, this will follow once we prove that
   $d(t)$ has a good factorization.

   If $d(t)$ is irreducible, then $d(t)=d(t)\cdot 1\cdot 1$ is a good
   factorization. For the rest of the proof we assume that $d(t)$ is
   reducible. If $d(t)$ has a linear factor, then $d(-1)=0$. Indeed,
   the only possible rational roots of $d(t)$ are $\pm 1$, and
   $d(1)=0$ is excluded by \ref{embedded}. For each type of polynomial
   in the chart, we study the factorization of $d(t)$ in the two
   remaining cases: $d(-1)=0$ and $d(t)$ has no linear factor.

  $\G(l+1)$ or $\GGO$: If $d(-1)=0$, then $l=4$ or $l=8$.  Thus, we
  have good factorizations $d(t)=(1-3t+t^2)\cdot(1+t)^2\cdot 1$ or
  $d(t)=(1-4t+t^2)\cdot (1+t)^2\cdot 1$\,.

  If $d(t)$ has no linear factor, then $d(t)=(1+at+\varepsilon
  t^2)(1+bt+\varepsilon t^2)$ with $a,b\in\mathbb Z$, both factors
  irreducible, and $\varepsilon=\pm 1$.  If $\varepsilon=-1$, then
  comparison of the coefficients of $t$ and $t^3$ gives $a+b<0$ and
  $-a-b<0$, so this case does not occur.  If $\varepsilon=1$, then
  comparison of the coefficients of $t^2$ gives $ab<0$, hence we may
  assume $b>0$, and then $d(t)=(1+at+t^2)\cdot(1+bt+t^2)\cdot 1$ is a
  good factorization.

 $\GTE$: If $d(-1)=0$, then $l=6$ and $d(t)=(1-4t+3t^2-t^3)\cdot(1+t)^2\cdot 1$
  is a good factorization.

  If $d(t)$ has no linear factor, then $d(t)=(1+at+\varepsilon
  t^2)(1+bt+ct^2-\varepsilon t^3)$ with $a,b,c\in \mathbb Z$, both
  factors irreducible, and $\varepsilon=\pm 1$. Comparing
  coefficients, we get:
\begin{align*}
a+b&=-2\\
ac+\varepsilon b-\varepsilon&=1\\
\varepsilon c-\varepsilon a&=1
\end{align*}
Using the first and the last equality to eliminate $b$ and $c$ from
the middle, we obtain $a^2=1+3\varepsilon$. If $\varepsilon=-1$, then
$a^2=-2$, which is not possible. If $\varepsilon=1$, then $a^2=4$,
hence $a=\pm 2$, and this contradicts the assumption that
$1+at+\varepsilon t^2$ is irreducible.

$\GH(p)$: If $d(-1)=0$, then $p+l=10$ and $d(t)=e(t)(1+t)$, where
$e(t)= 1-3t+(p-5)t^2+3t^3-t^4$. Since $1+t$ has non-negative
coefficients, the polynomial $d(t)$ has a good factorization if and
only if $e(t)$ has a good factorization. If $e(t)$ is irreducible,
then $e(t)=e(t)\cdot 1\cdot 1$ is a good factorization, so we assume
$e(t)$ is reducible.  If $e(t)$ has a linear factor, then $e(-1)=0$,
hence $p=5$. It follows that $p=l=5$, and this is ruled out by
\ref{embedded}. If $e(t)$ has no linear factor, then
$e(t)=(1+at-t^2)(1+bt+t^2)$ with both factors irreducible and $a,b\in
\mathbb Z$. Comparing the coefficients of $t$ and $t^3$, we get
$a+b=-3$ and $a-b=3$, hence $a=0$. This contradicts the hypothesis
that $1+at-t^2$ is irreducible.

If $d(t)$ has no linear factor, then $d(t)=(1+at+\varepsilon
t^2)(1+bt+ct^2-\varepsilon t^3)$ with $a,b,c\in \mathbb Z$, both
factors irreducible, and $\varepsilon=\pm 1$. Comparison of
coefficients yields:
\begin{align*}
a+b&=-2\\
ab+c+\varepsilon&=-l+2\\
ac+\varepsilon b-\varepsilon&=p-2\\
\varepsilon c-\varepsilon a&=2
\end{align*}
These equalities yield $-a^2-a=-l+2-3\varepsilon$ and $a^2+\varepsilon
a=p-2+3\varepsilon$, so $\varepsilon a-a=p-l$.  By \ref{classify} and
\ref{embedded} we have $p<l$, hence $\varepsilon=-1$ and $2a=l-p>0$.
Thus, $d(t)=(1+bt+ct^2+t^3)\cdot 1\cdot(1+at-t^2)$ is a good
factorization.
 \end{proof}

\section{Finite CI-dimension}

In this section we let $R$ denote a local ring and let $M$, $N$ be
finite $R$-modules.

\begin{chunk}
We refer to the introduction for the definition of finite
CI-dimension and we recall below the basic examples:

\begin{subchunk}
\mlabel{ex2} If $R$ is a complete intersection, then $M$ has
finite CI-dimension.
\end{subchunk}

\begin{subchunk}
\mlabel{ex1} If $\pd_RM<\infty$, then $M$ has finite CI-dimension.
\end{subchunk}
\end{chunk}

\begin{chunk}
\mlabel{rmk} When one of the modules $M$ or $N$ has finite
CI-dimension, several (in)equalities involving the numbers
$e_R(M,N)$ and $t^R(M,N)$ are known.

The inequality below follows from \cite[4.9]{AB} and
\cite[1.4]{AGP}:

\begin{subchunk}
\mlabel{rmk1} If $t^R(M,N)<\infty$ and $M$ has finite
CI-dimension, then
\[
t^R(M,N)\le\depth R-\depth_RM\,.
\]
\end{subchunk}

A formula for $e_R(M,N)$ is given by Araya and Yoshino
\cite[4.2]{AY}:

\begin{subchunk}
\mlabel{rmk2} If $e_R(M,N)<\infty$ and $M$ has finite
CI-dimension, then
\[
e_R(M,N)=\depth R-\depth_RM\,.
\]
\end{subchunk}
\end{chunk}
Over Gorenstein rings the equality of \ref{rmk2} is valid more
generally:

\begin{lemma}
\mlabel{rmk3} If $R$ is Gorenstein, $e_R(M,N)$ is finite and $N$
has finite CI-dimension, then
\[
e_R(M,N)=\depth R-\depth_RM\,.
\]
\end{lemma}

\noindent{\em Proof.\/} If $\pd_RN<\infty$, then $\id_RN<\infty$
(see \ref{Gor}), so the equality is given by a result of Ischebeck
\cite[2.6]{I}. In general, there exists a flat ring homomorphism
$R\to R'$ and a deformation $Q$ of $R'$ such that
$\pd_Q(M\otimes_RR')<\infty$. By \ref{completion} we may assume
that $R=R'$ and then $e_Q(M,N)=\depth Q-\depth_QM$. A standard
argument, cf.\ \cite[2.6]{AY}, then gives
\begin{xxalignat}{3}
&\ &e_R(M,N)&=e_Q(M,N)-(\depth_QR-\depth Q)=\depth
R-\depth_RM\,.&&\square
\end{xxalignat}

The next theorem is the main result of this section. It shows, in
particular, that every Gorenstein ring of codimension at most $4$
is AB.

\begin{theorem}
\mlabel{T2} Let $R$ be a local Gorenstein ring with $\codim R\le
4$ and let $M$, $N$ be finite $R$-modules.

If $e_R(M,N)$ is finite, then $M$ or $N$ has finite CI-dimension
and
\begin{align*}
e_R(M,N)&=\dim R-\depth_RM\\
e_R(N,M)&=\dim R-\depth_RN
\end{align*}

If $t^R(M,N)$ is finite, then $M$ or $N$ has finite CI-dimension
and
\[
t^R(M,N)\le\max\{\dim R-\depth_RM,\dim R-\depth_RN\}\,.
\]
\end{theorem}

\begin{proof}

In view of the results recalled in \ref{rmk}, Lemma \ref{rmk3},
and Theorem HJ stated in the introduction, it suffices to prove
the following claim: If $e_R(M,N)$ or $t^R(M,N)$ is finite, then
one of the modules $M$, $N$ has finite CI-dimension.

By \ref{ex2}, the claim holds when $R$ is a complete intersection.
If $R$ has no embedded deformation, then the statement results
from Theorem \ref{T1}, in view of \ref{ex1}.

It remains thus to treat the case when $R$ is not complete
intersection and has an embedded deformation $Q$. This only
happens when $\codim R=4$ and $Q$ has no embedded deformation
(otherwise, $\wh R$ deforms to a Gorenstein ring $Q'$ with $\codim
Q'\le 2$, and such a ring is a complete intersection).  Standard
arguments (see \cite[2.6]{AY} for example) show that
$e_R(M,N)<\infty$ implies $e_Q(\wh M,\wh N)<\infty$ and $t^R(
M,N)<\infty$ implies $t^Q(\wh M,\wh N)<\infty$. By Theorem 1, $\wh
M$ or $\wh N$ has then finite projective dimension over $Q$, and
thus the corresponding module has finite CI-dimension over $R$.
\end{proof}

We note some further applications of the conclusion on finite
CI-dimension of Theorem \ref{T2}.

Araya and Yoshino \cite[4.2]{AY} give a self-test for finite
projective dimension when the CI-dimension is finite: If $M$ has
finite CI-dimension, then there is an equality
$e_R(M,M)=\pd_R(M)$. In view of Theorem \ref{T2} we have thus:

\begin{corollary}
There is an equality $e_R(M,M)=\pd_R(M)$. \qed
\end{corollary}

When $M$ or $N$ has finite CI-dimension, (1) below is proved by
Jorgensen \cite[2.2]{J} and (2) is proved by Araya and Yoshino
\cite[2.5]{AY} (see also Iyengar \cite[4.3]{Iy} for the case
$q=0$) over any local ring $R$. By Theorem \ref{T2} we have thus:

\begin{corollary}
Set $d_R(M,N)=\depth R-\depth_RM-\depth_RN$.

 If
$t^R(M,N)=q<\infty$, then the following hold:
\begin{enumerate}[{\quad\rm(1)}]
\item
$t^R(M,N)=\sup\{d_{R_\fp}(M_\fp,N_\fp)\mid\fp\in\supp M\cap\supp
N\}$.
\item If $q=0$ or $\depth_R\Tor^R_{q}(M,N)\le 1$,
 then
\begin{xxalignat}{3}
&\ &t^R(M,N)&=d_R(M,N)+\depth_R\Tor_q^R(M,N)\,.&&\square
\end{xxalignat}
\end{enumerate}
\end{corollary}

The conclusion on CI-dimension of Theorem \ref{T2} does not extend
to higher codimensions: There exist rings of any codimension
greater than or equal to $6$ and $R$-modules $M$, $N$ of infinite
CI-dimension such that $e_R(M,N)<\infty$. The example below is
based on a construction in \cite[4.3]{HJ}.

\begin{example}
Let $(S,\fs,k)$ and $(T,\ft,k)$ be two Gorenstein rings that are
essentially of finite type over a field $k$ and none of them is a
complete intersection, so they have codimension at least $3$. Set
$\fp=\fs\otimes_kT+S\otimes_k\ft$. The local ring
$R=(S\otimes_kT)_{\fp}$ is then Gorenstein and has $\codim
R=\codim T+\codim S\ge 6$. The $R$-modules $M=S\otimes_kk$ and
$N=k\otimes_kT$ satisfy $e_R(M,N)\le\dim R$ by \cite[4.3]{HJ}.
Since $S$ is not a complete intersection, $k$ has infinite
CI-dimension over it, cf.\ \cite[1.3]{AGP}. The ring $S$ is
faithfully flat as a $k$-module, so \cite[1.13]{AGP} yields
$\CIdim_RM\ge\CIdim_Sk=\infty$. By symmetry, $N$ also has infinite
CI-dimension over $R$.
\end{example}

 \section{Embedded deformations}

 In this section we construct examples to show that vanishing of
 (co)homology can occur for nontrivial reasons.

 Throughout the section, $S$ denotes a local commutative noetherian
 ring. Our examples are based on the existence of modules with
 periodic resolutions over rings with embedded deformations. We recall
 the relevant definitions and results:

 Let $a\ge 1$ be an integer. The minimal free resolution of a finite
 $S$-module $U$ is said to be {\em periodic of period $a$} if
 $\Syz_n^S(U)\cong\Syz_{n+a}^S(U)$ for all $n\ge 0$. The following
 result is due to Avramov, Gasharov and Peeva \cite[3.2]{AGP}:

 \begin{chunk}
   \mlabel{periodic} If $S\cong Q/(f)$ for a non-regular local ring $(Q,\fn)$ and a
   $Q$-regular element $f\in\fn^2$, then there exists a finite
   $S$-module $U$ whose minimal $S$-free resolution is periodic of
   period $2$ and $\pd_Q U=1$. In particular, $\pd_S U=\infty$ and
   $\depth_S U=\depth S$.
 \end{chunk}

 \begin{theorem}
   \mlabel{deformed} If $S\cong Q/(f)$ for a non-regular local ring $(Q,\fn)$ and a
   $Q$-regular element $f\in\fn^2$, and $U$ is the $S$-module from
   {\rm\ref{periodic}}, then there exists a finite $S$-module $V$ with
   $\depth_S V=\depth S$ such that
\[
   \text{$\id_SV=\infty=\pd_S V$\quad and \quad$e_S(U,V)=0=t_S(U,V)$.}
\]

If, furthermore, $S$ is Gorenstein, then also $e_S(V,U)=0$.
\end{theorem}

 \begin{proof}
   Since $Q$ is not regular, we have $\pd_Qk=\infty$. Set $s=\depth
   Q+2$. By \cite[\S2]{O} we have $\depth_Q\Syz^Q_i(k)=\depth Q$ for
   each $i\ge s$. If $Q$ is Gorenstein, then $\id_Q\Syz^Q_i(k)=\infty$
   for each $i\ge 0$ by \ref{Gor} and we set $V'=\Syz^Q_s(k)$.
   Otherwise, one of the syzygy modules in the short exact sequence
   $0\to\Syz^Q_{s+1}(k)\to Q^{b}\to\Syz^Q_{s}(k)\to 0$ has infinite
   injective dimension and we let $V'$ be this module.

   Set $V=V'/fV'$. Note that $f$ is a $V'$-regular element and
   $\depth_S V=\depth S$. The isomorphisms
   $\Ext^n_{S}(k,V)\cong\Ext^{n+1}_Q(k,V')$ and
   $\Tor^{S}_n(k,V)\cong\Tor^Q_n(k,V')$ show that
   $\id_{S}V=\pd_{S}V=\infty$. The isomorphisms
   $\Ext^n_{S}(U,V)\cong\Ext^{n+1}_Q(U,V')$ and
   $\Tor^{S}_n(U,V)\cong\Tor^Q_n(U,V')$ yield $e_{S}(U,V)=0$ and
   $t^{S}(U,V)\le 1$, because $\pd_QU=1$. The periodicity of the
   minimal free resolution of $U$ gives
   $\Tor^{S}_n(U,V)\cong\Tor^{S}_{n+2}(U,V)$ for all $n>0$, hence
   $t^{S}(U,V)=0$.

   If $S$ is a Gorenstein ring, then $Q$ is Gorenstein, hence
   $\id_QU<\infty$ by \ref{Gor}. The isomorphisms
   $\Ext^n_{S}(V,U)\cong\Ext^{n}_Q(V',U)$ then show
   $e_{S}(V,U)<\infty$. Since the $S$-module $U$ has finite CI-dimension, Lemma \ref{rmk3} yields
   $e_S(V,U)=0$.
 \end{proof}

 The next lemma is an extension of the fact that if $U$ is a finite
 $S$-module and $\ov U=U/gU$ for an $U$-regular element $g$, then
 $\pd_S \ov U=\pd_S U+1$ and $\id_S \ov U=\id_S U$.

 \begin{lemma}
   \mlabel{artinian} Let $S$ be a local ring, $U$, $V$ be finite
   $S$-modules and $g\in S$ an $U$-regular element. For the $S$-module
   $\ov U=U/gU$ the following hold:
 \begin{enumerate}[{\quad\rm(1)}]
 \item $e_S(\ov U,V)=e_S(U,V)+1$.
 \item $e_S(V,\ov U)=e_S(V,U)$.
 \item $t^S(U,V)\le t^S(\ov U,V)\le t^S(U,V)+1$.\\
   If $gV=0$ or $V$ has finite length, then $t^S(\ov U,V)=t^S(U,V)+1$.
 \item If the minimal free resolution of $\Syz_s^S(U)$ is periodic of
   period $2$ for an integer $s\ge 0$, then the minimal free
   resolution of $\Syz_{s+1}^S(\ov U)$ is periodic of period $2$.
\end{enumerate}
 \end{lemma}

 \begin{proof}
   Properties (1) to (3) are deduced using the long exact sequences
   induced by the exact sequence $0\to U\xrightarrow{g} U\to \ov U\to
   0$.

   We only give the proof of (3). In this case, the long exact
   sequence is:
 \begin{align*}
   \hdots\to\Tor_{n+1}^S(\ov U,V)&\to\Tor_{n}^S(U,V)\xrightarrow{g}\Tor_{n}^S(U,V)\to\\
   &\to\Tor_{n}^S(\ov U,V)\to\Tor_{n-1}^S(U,V)\to\hdots
 \end{align*}
 Nakayama's Lemma shows that if $\Tor_{n}^S(U,V)\ne 0$ then
 $\Tor_{n}^S(\ov U,V)\ne 0$, hence $t^S(U,V)\le t^S(\ov U,V)$. To
 prove the remaining statements it suffices to assume
 $t^S(U,V)=p<\infty$. The long exact sequence shows that if
 $\Tor^S_i(U,V)=0$ for $i=n$ and $i=n-1$ then $\Tor_{n}^S(\ov U,V)=0$.
 We conclude $t^S(\ov U,V)\le p+1$. If $gV=0$ or $V$ has finite
 length, then multiplication by $g$ on $\Tor_{p}^S(U,V)$ has non-zero
 kernel, hence $\Tor_{p+1}^S(\ov U,V)\ne 0$ and thus $t^S(\ov
 U,V)=p+1$.

 (4) Let $F$ be a minimal free resolution of $U$. The mapping cone of
 the homomorphism $F\xrightarrow{g} F$ is a minimal free resolution of
 $\ov U$, hence we have $\Syz_{n}^S(\ov
 U)\cong\Syz_{n}^S(U)/g\Syz_{n}^S(U)\oplus\Syz_{n-1}^S(U)$ for all
 $n\ge 0$, and the conclusion follows.
\end{proof}

 \begin{corollary}
   \mlabel{last} Let $R$ be a $d$-dimensional Cohen-Macaulay local
   ring such that $\codim R\ge 2$ and $\wh R\cong Q/(f)$ for a local
   ring $(Q,\fn)$ and a $Q$-regular element $f\in\fn^2$.

   There exist $R$-modules
   $M$, $N$ of finite length such that $\pd_Q M=d+1$, the minimal
   $R$-free resolution of $\Syz^R_{d}(M)$ is periodic of period $2$ and
\[
\text{$\id_R N=\infty=\pd_R N$ \quad and\quad $e_R(M,N)=d=t^R(M,N)$.}
\]

If, furthermore, $R$ is Gorenstein, then $e_R(N,M)=d$.

 \end{corollary}

 \begin{proof}
   Since $\codim R\ge 2$, the ring $Q$ is not regular. Set $S=\wh R$
   and let $U$ and $V$ be as in Theorem \ref{deformed}. Choose a $U$-
   and $V$-regular sequence $\bd g$ of length $d$ and set $M=U/\bd gU$
   and $N=V/\bd gV$. By Lemma \ref{artinian} we have $\id_S
   N=\infty=\pd_S N$ and the minimal $S$-free resolution of
   $\Syz^S_{d}(M)$ is periodic of period $2$.  Using Theorem
   \ref{deformed} and Lemma \ref{artinian} we conclude
   $e_{S}(M,N)=e_S(U,V)+d=d$ and $t^{S}(M,N)=t^{S}(U,N)+d<\infty$.
   Since the minimal free resolution of $U$ is periodic of period $2$,
   we have $\Tor^{S}_n(U,N)\cong\Tor^{S}_{n+2}(U,N)$ for all $n>0$,
   hence $t^{S}(U,N)=0$ and thus $t^{S}(M,N)=d$. If $S$ is a
   Gorenstein ring, then $e_{S}(N,M)=e_{S}(V,U)+d=d$ by Theorem
   \ref{deformed} and Lemma \ref{artinian}.

   It remains to notice that $M$ and $N$ are $\wh R$-modules of finite
   length, hence, considered as $R$-modules, they have finite length
   and are complete. The desired conclusion now follows by applying
   \ref{completion}
 \end{proof}

\section*{Acknowledgment}
I wish to thank my thesis advisor Luchezar Avramov for inspiring the
development of this paper and for his constant help.

\end{document}